\documentclass[reqno]{amsart}
\pdfoutput=1
\usepackage{amssymb}
\usepackage{graphicx}
\usepackage[all]{xy}

\begin{document}
\newcommand*{\threesim}{\mathrel{\vcenter{\offinterlineskip\hbox{$\sim$}\vskip-.35ex\hbox{$\sim$}\vskip-.35ex\hbox{$\sim$}}}}
\newtheorem{prop}{Proposition} [section]
\newtheorem{lemma}[prop]{Lemma}
\newtheorem{thm}[prop]{Theorem}
\newtheorem{cor}[prop]{Corollary}
\newtheorem*{thmA}{Theorem A}
\newtheorem*{thmB}{Theorem B}
\newtheorem*{thmC}{Theorem C}
\theoremstyle{definition}
\newtheorem{exa}{Example}
\newtheorem*{rmk}{Remark}
\theoremstyle{remark}
\newtheorem*{ack}{Acknowledgment}
\newtheorem*{pf}{Proof}
\numberwithin{equation}{section}
\title{The Lind Zeta functions of reversal systems of finite order}
\author{Sieye Ryu}
\address{Department of Mathematics Ben Gurion University of the Negev,
P.O.B. 653, Be'er Sheva, 8410501, Israel}
\email{Sieye Ryu <sieye@math.bgu.ac.il>}

\maketitle
\begin{abstract}
A decomposition theorem for the Lind zeta function of a reversal system $(X, T, R)$ of finite order is established. A reversal system can be regarded as an action of a certain group $G$ on $X$. To establish an explicit formula for the Lind zeta function of $(X, T, R)$, we need to consider finite index subgroups $H$ of $G$ with induced actions given by automorphisms or by flips. 

When the underlying dynamical system $(X, T)$ is either a shift of finite type or a sofic shift, we express the Lind zeta function of $(X, T, R)$ in terms of matrices.  
\end{abstract}

\section{Introduction} 
The concept of time reversibility plays an important role in mathematics, physics and computer science.
For instance, in celestial mechanics, time reversibility has been studied in order to investigate the initial conditions and history of the universe. In thermodynamics, time reversibility has been studied due to its connections with energy efficiency and entropy. The application of entropy and time reversibility has also attracted attention in information theory. 
Time reversibility is a natural feature of many areas of mathematics, including dynamical systems \cite{LR, Bi, 1S}, ergodic theory \cite{G}, and automaton theory \cite{Ka}.

In \cite{LPS}, Lee et all studied the notion of reversible shift dynamical systems. 
A shift space is a shift-invariant and closed subset of $\mathcal{A}^{\mathbb{Z}}$, where $\mathcal{A}$ is a finite set. Here, we assume that  $\mathcal{A}$ is given the discrete topology and $\mathcal{A}^{\mathbb{Z}}$ is given the product topology. 
A shift space with the shift map can be regarded as an action of $\mathbb{Z}$. In order to generalize this principle, one can also study the actions of other groups on a shift space. For instance, actions of $\mathbb{Z}^d$ and other groups on shift spaces have been studied in \cite{HM, Ki, L} and \cite{AKM, BS, F, 0KR1}, respectively. Time reversal for a shift space can be regarded as a non-abelian group action.
In this article, we study the Lind zeta function of a reversal system and express it in terms of matrices when its underlying dynamical system is either a shift of finite type or a sofic shift.
For the definitions of a shift of finite type and a sofic shift, see \cite{LM}.

We begin by defining reversal systems.     
Throughout the article, we assume that $(X, T)$ is an invertible dynamical system. 
A homeomorphism $R:X \rightarrow X$ is said to be a \textit{reversal for $(X, T)$} if
\begin{equation} \label{eq: 1.1}
T \circ R = R \circ T^{-1}.
\end{equation}
The triple $(X, T, R)$ is called a \textit{reversal system}. 

Suppose that $(X, T, R)$ is a reversal system. 
We call the number $|R|=\min\{n\in \mathbb{N}\setminus\{0\} : R^n = \text{id}_X\}$ the \textit{order} of $R$.
If $|R|=r$, then $(X, T, R)$ is said to be a \textit{reversal system of order $r$}.
If $R$ is of order $r$ for some positive integer $r$, 
then  $R$ is said to be a \textit{reversal of finite order}
and the triple $(X, T, R)$ is said to be a \textit{reversal system of finite order}. We will see in Section 2, if $|R|$ is finite, then it must be an even number unless $T=T^{-1}$. From now on, we assume that the order of our reversal is $2r$.  

Suppose that $G_{2r}$ is a group with presentation $\langle a, b : ab=ba^{-1} \, \text{ and } \, b^{2r}=1 \rangle$.
To each reversal system $(X, T, R)$ of order $2r$, there corresponds a unique $G_{2r}$-action $\alpha_{2r}$ on $X$ such that 
$$\alpha_{2r}(a, x) = T(x) \qquad \text{and} \qquad \alpha_{2r}(b, x) = R(x).$$
Two reversal systems $(X, T, R)$ and $(X', T', R')$ are said to be $G_{2r}$-\textit{conjugate} if there is an equivariant homeomorphism $\theta: (X, T, R) \rightarrow (X', T', R'),$ that is, $\theta$ satisfies the following:
\begin{equation}\label{eq: 1.2}
\theta \circ T = T' \circ \theta \qquad \text{and} \qquad \theta \circ R = R' \circ \theta
\end{equation}
and we write $(X, T, R) \cong (X', T', R')$. The homeomorphism $\theta$ is called a $G_{2r}$-\textit{conjugacy from $(X, T, R)$ to $(X', T', R')$}.
It is obvious that two $G_{2r}$-conjugate reversal systems have the same order.
 
Suppose that $G$ is a group, $X$ is a set and $\alpha : G \curvearrowright X$ is a $G$-action on $X$.  
Let $\mathcal{F}$ denote the set of finite-index subgroups.
For each $H \in \mathcal{F}$, we denote the number of points fixed by all elements $h\in H$ by $f_{\alpha}(H)$:
$$f_{\alpha}(H) = |\{ x \in X : \forall h \in H \; \alpha(h, x) = x  \}|.$$
The Lind zeta function \cite{L} is formally defined by
\begin{equation}\label{eq: 1.3}
\zeta_{\alpha}(t) = \exp \left( \sum_{H \in \mathcal{F}} \;\frac{f_{\alpha}(H)}{|G/{H}|}t^{|G/H|}\right).
\end{equation}
Suppose that $(X, T, R)$ and $(X', T', R')$ are reversal systems of order $2r$. Let $\alpha:G_{2r} \curvearrowright X$ and $\alpha':G_{2r}\curvearrowright X'$ be the corresponding $G_{2r}$-actions on $X$ and $X'$, respectively. It is obvious that $(X, T, R) \cong (X', T', R')$ implies $f_{\alpha}(H)=f_{\alpha'}(H)$. As a consequence, the Lind zeta function is a $G_{2r}$-conjugacy invariant.

The first goal of this article is to establish an explicit formula for the Lind zeta function of a reversal system of order $2r$.
In \cite{KLP}, Kim et all investigated the case when $r=1$. In this article, we will prove the following theorem.

\begin{thmA}
The Lind zeta function $\zeta_{T, R}$ of a reversal system $(X, T, R)$ of order $2r$ is given by
$$\zeta_{T, R}(t) = \prod_{\substack{k|r \\ 1 \leq k \leq r}}\, \exp \Big(\frac{g_{2k}(t^{2k})}{2k} \Big)\prod_{\substack{2k-1|r \\ 1 \leq 2k-1 \leq r} }\exp \left( \frac{h_{4k-2}(t^{2k-1})}{2k-1} \right).$$
Here, $g_{2k}(t)$ are the generating functions of the sub-reversal systems of $(X, T, R)$, which have $R$ as reversals of order $2k$
and $h_{4k-2}(t)$ are the generating functions of the sub-reversal systems of $(X, T, R)$, which have $R^{2k-1}$ as reversals of order $2$.
\end{thmA}

We briefly describe what those generating functions are and will give the precise definitions in Section 3. 
As we will see in Section 3, if $H$ is a finite index subgroup of $G_{2r}$, then the induced actions is either by automorphisms or by reversals.
The generating functions $g_{2k}(t)$ record the numbers of fixed points of the actions by automorphisms. The generating functions $h_{4k-2}(t)$ record the same information when the actions are by reversals of order $2$. 

The second goal of this article is to express the Lind zeta function of a  reversal system $(X, \sigma_X, \varphi)$ in terms of matrices when $(X, \sigma_X)$ is either a shift of finite type or a sofic shift. 
As we will see in Section 3, Theorem A implies that it suffices to consider the numbers of fixed points of actions $G \curvearrowright X$ either by automorphisms or by reversals of order $2$.
The case when actions are by reversals of order $2$ is investigated in \cite{KLP, KR}. (See also Section 2.) In this article, we focus on the case when actions are by automorphisms. 

We introduce some notation.
A reversal system $(X, \sigma_X, \varphi)$ will be called a \textit{shift-reversal system of finite type} if $(X, \sigma_X)$ is a shift of finite type and it will be called a \textit{sofic shift-reversal system} if $(X, \sigma_X)$ is a sofic shift.

If $m$ and $l$ are positive integers with $0 \leq l <r$, then
$f_{\sigma_X, \varphi}(m, 2l)$ will denote the number of points in $X$ fixed by $(\sigma_X)^m \circ {\varphi}^{2l}$:
$$
f_{\sigma_X, \varphi}(m, 2l) = |\{x \in X : (\sigma_X)^m \circ {\varphi}^{2l}(x) = x\}|.
$$

The second goal will be achieved by Theorem B and Theorem C.  

\begin{thmB}
Suppose that $(X, \sigma_X, \varphi)$ is a shift-reversal system of finite type. Then there are zero-one square matrices $A$ and $J$ such that
$$f_{\sigma_X, \varphi}(m, 2l) =\emph{tr}(A^m J^{2l}) \qquad (m=1,2, \cdots, \; 0 \leq l <r).$$
\end{thmB} 

\begin{thmC}
If $(X, \sigma_X , \varphi)$ is a sofic shift-reversal system, then there are square matrices $A_k$ and $J_k, k = 1, 2, . . . , n$, whose entries are integers such that
$$f_{\sigma_X, \varphi}(m, 2l) =\sum_{k=1}^{n}\, (-1)^{k+1}\emph{tr}((A_k)^m (J_k)^{2l}) \qquad (m=1,2, \cdots, \; 0 \leq l <r).$$
\end{thmC}

As a consequence of Theorem B and Theorem C, we obtain the following corollary.

\begin{cor}
Suppose that $(X, \sigma_X, \varphi)$ is either a shift-reversal system of finite type or a sofic-shift reversal system.
Then the generating functions $g_{2k}(t^{2k})$ are rational functions.
\end{cor}

In Section 6, we will see that the generating function $g_{2k}(t^{2k})$ is actually $\mathbb{N}$-rational. (See \cite{R, BR1, E} for the definition of an $\mathbb{N}$-rational function.) 
In \cite{KR}, it is proved that $h_{4k-2}(t^{2k-1})$ is $\mathbb{N}$-rational. However,
neither $\frac{g_{2k}(t^{2k})}{2k}$ nor $\frac{h_{4k-2}(t^{2k-1})}{2k-1}$
is $\mathbb{N}$-rational in general. 
We will discuss it in Section 6.

This article is organized as follows. In Section 2, we explain the background. We prove Theorem A, B, and C in Section 3, 4 and 5, respectively. In Section 6, there will be some remarks and an example.

\begin{ack}
{\rm The author obtained some partial results of this research during the author's doctoral course under the supervision of Young-One Kim. The author deeply thanks him for his patience. The author also thanks Jungseob Lee for his helpful comments. This research was partially supported by the Israel Science Foundation (grant no. 626/14) and the People Programme (Marie Curie Actions) of the European Union's Seventh Framework Programme (FP7/2007-2013) under REA grant agreement no. 333598.}
\end{ack}

\section{Preliminary} 
Suppose that $R$ is a reversal for $(X, T)$. 
By (\ref{eq: 1.1}), $R^{2k}$ is an automorphism of $(X, T)$ for all positive integers $k$ , that is, 
$T \circ R^{2k} = R^{2k} \circ T$.
On the other hand, $R^{2k-1}$ is a reversal for $(X, T)$ for all positive integers $k$.
Hence, if $|R|=r$ and $T^2 \neq \text{id}_X$, then $r$ must be even. 
From now on, $R$ will be said to be a \textit{flip} 
and the triple $(X, T, R)$ will be said to be a \textit{flip system} 
when $|R|=2$. When we say $(X, T, R)$ is a reversal system or $R$ is a reversal, we assume that $|R| \neq 2$. 

As in the case of reversals, a flip system $(X, \sigma_X, \varphi)$ will be called a \textit{shift-flip system of finite type} if $(X, \sigma_X)$ is a shift of finite type and it will be called a \textit{sofic shift-flip system} if $(X, \sigma_X)$ is sofic.

If $(X, T, R)$ is a reversal system, then so is $(X, T, T^n\circ R)$ for any integer $n$ by (\ref{eq: 1.1}). From (\ref{eq: 1.2}), it follows that $T^n$ is a conjugacy from $(X, T, R)$ to $(X, T, T^{2n}\circ R)$.

Lind introduced the function in (\ref{eq: 1.3}) in the special case where $G=\mathbb{Z}^d$ in \cite{L} and it is a generalization of the Artin-Mazur zeta function, which is introduced in \cite{AM}. We recall the Artin-Mazur zeta function.

If $m$ is a positive integer, then the number of periodic points in $(X, T)$ of period $m$ is denoted by $p_m(T)$:
$$p_T(m) = | \{ x \in X : T^m(x)=x \}|.$$
Suppose that the sequence $\{(p_T (m))^{1/m}\} $ is bounded. Then the Artin-Mazur zeta function $\zeta_{T}$ is defined by
$$\zeta_{T}(t) = \exp \biggl( \sum_{m=1}^{\infty} \frac{p_T(m)}{m} t^m \biggr).$$ 
One can prove that if $\alpha: \mathbb{Z} \curvearrowright X$ is given by
$(m, x) \mapsto T^mx$, then the Lind zeta function $\zeta_{\alpha}$ becomes the Artin-Mazur zeta function $\zeta_T$.

In \cite{KLP}, an explicit formula for the Lind zeta function $\zeta_{T, F}$ of a flip system $(X, T, F)$ is established. 
If $m$ is a positive integer and $n$ is an integer, then
$p_{T,F} (m, n)$ will denote the number of points in $X$ fixed by $T^m$ and $T^n\circ F$:
$$p_{T, F} (m, n) = |\{x \in X : T^m (x) = T^n \circ F(x) = x\}|.$$
The Lind zeta function of $(X, T, F)$ is given by
\begin{equation}\label{eq: 2.1}
\zeta_{T, F}(t) = \sqrt{\zeta_{T}(t)} \exp \left(h_{T, F}(t) \right),
\end{equation}
where $\zeta_{T}$ is the Artin-Mazur zeta function  
and
\begin{equation}\label{eq: 2.11}
h_{T, F}(t) = \sum_{m=1}^{\infty} \,\left(p_{T, F}(2m-1, 0)t^{2m-1} + \frac{p_{T, F}(2m, 0)+p_{T, F}(2m, 1)}{2} t^{2m}\right).
\end{equation}
The function $h_{T, F}(t)$ is called the generating function.
In Section 3, we will sketch the proof of (\ref{eq: 2.1}).

It is well known \cite{LM} that if $X$ is a shift of finite type, then there is a square matrix $A$ whose entries are non-negative integers such that
\begin{equation}\label{eq: 2.2}
p_{\sigma_X}(m) = \text{tr}(A^m )\qquad (m = 1, 2, \cdots).
\end{equation}
It is also well known \cite{B, M, LM} that if $X$ is a sofic shift, then there are square matrices $A_1 , A_2 , . . . , A_n$ whose entries are integers such that
\begin{equation}\label{eq: 2.3}
p_{\sigma_X}(m) = \sum_{k-1}^n \, (-1)^{k+1}\, \text{tr}((A_k)^m ) \qquad
(m = 1, 2, \cdots).
\end{equation}
Equations (\ref{eq: 2.2}) and (\ref{eq: 2.3}) allow us to express the Artin-Mazur zeta function $\zeta_{\sigma_X}$ of $(X, \sigma_X)$ in terms of matrices when $(X, \sigma_X)$ is either a shift of finite type or a sofic shift. If $X$ is a shift of finite type, then $\zeta_{\sigma_X}$ is the reciprocal of a polynomial and if $X$ is sofic, then $\zeta_{\sigma_X}$ is a rational function. In either cases,  $\zeta_{\sigma_X}$ is $\mathbb{N}$-rational \cite{R, BR1, E}.

Similarly, when $(X, \sigma_X, \varphi)$ is either a shift-flip system of finite type or a sofic shift-flip system, the number of fixed points $p_{\sigma_X, \varphi}(m, n)$ can be expressed in terms of matrices \cite{KLP, KR}.
To present it, we need notation.
If $M$ is a matrix, $S[M]$ will denote the sum of the entries of $M$, that is,
$$\mathcal{S}[M] = \sum_{I,J} M_{I J}$$
and $M^{\triangle}$ will denote the matrix whose diagonal entries are identical with those of $M$ but whose other entries are equal to $0$, that is,
$$(M^{\triangle})_{IJ} =  \begin{cases} M_{IJ} \qquad \text{if } I=J, \\ 0 \qquad \quad \text{otherwise}. \end{cases} $$

The following is proved in \cite{KLP}.
If $(X, \sigma_X , \varphi)$ is a shift-flip system of finite type, then there are zero-one square matrices $A$ and $J$ such that
\begin{eqnarray*}
p_{\sigma_X , \varphi}(2m-1, 0) &=& \mathcal{S}[J^{\triangle} A^{m-1} (AJ)^{\triangle} ], \\
p_{\sigma_X , \varphi}(2m, 0) &=& \mathcal{S}[J^{\triangle} A^{m} J^{\triangle} ] \qquad \text{and} \\
p_{\sigma_X , \varphi}(2m, 1) &=& \mathcal{S}[(JA)^{\triangle} A^{m-1} (AJ)^{\triangle} ] \qquad
(m = 1, 2, \cdots). 
\end{eqnarray*}

The following is proved in \cite{KR}.
If $(X, \sigma_X , \varphi)$ is sofic shift-flip system, then there are square matrices $A_k, B_k$ and $J_k, k = 1, 2, . . . , n$, whose entries are integers such that
\begin{eqnarray}\label{eq: 2.5}
p_{\sigma_X , \varphi}(2m-1, 0) &=& \sum_{k=1}^n \, (-1)^{k+1} \mathcal{S}[J_k^{\triangle} A_k^{m-1} (A_kJ_k)^{\triangle} ], \nonumber \\
p_{\sigma_X , \varphi}(2m, 0) &=& \sum_{k=1}^n \, (-1)^{k+1} \mathcal{S}[J_k^{\triangle} A_k^{m} J_k^{\triangle} ] \qquad \text{and} \\
p_{\sigma_X , \varphi}(2m, 1) &=& \sum_{k=1}^n \, (-1)^{k+1} \mathcal{S}[(J_kA_k)^{\triangle} A_k^{m-1} (A_kJ_k)^{\triangle} ] \qquad
(m = 1, 2, \cdots).\nonumber
\end{eqnarray}

From the above formulas, 
one can see that the generating functions of a shift-flip system of finite type and a sofic shift flip systems are rational functions. In \cite{KR}, it is shown that they are actually $\mathbb{N}$-rational.

\section{Proof of Theorem A}

We first sketch the proof of (\ref{eq: 2.1}). For more details, see \cite{KLP}. 
Every finite index subgroup of $G_2 =\langle a, b : ab=ba^{-1} \;\, \text{and} \;\, b^2=1\rangle$ can be written in one and only one of the following forms: 
$$\langle a^m \rangle \qquad \text{or} \qquad \langle a^m, a^k b \rangle \qquad (m=1, 2, \cdots , \; k=1, 2, \cdots, m-1)$$
and has index 
$$|G_2/ \langle a^m \rangle| = 2m \qquad \text{or} \qquad |G_2/ \langle a^m, a^k b \rangle| = m.$$
From this, we obtain
$$\zeta_{T, F}(t) = \exp \Big(\sum_{m=1}^{\infty} \, \frac{p_T(m)}{2m}t^{2m} +\sum_{m=1}^{\infty} \, \sum_{k=0}^{m-1}\, \frac{p_{T,F}(m,k)}{m}t^{m} \Big).$$
The definition of a flip tells us that
$$p_{T, F}(m,n) = p_{T, F}(m,n+m) = p_{T, F}(m,n+2)$$
and this implies 
$$p_{T, F}(m,n) = \begin{cases} \, p_{T, F}(m, 0) \qquad \text{ if } m \text{ is odd}, \\ \, p_{T, F}(m, 0) \qquad \text{ if } m \text{ and } n \text{ are even}, \\ \, p_{T, F}(m, 1) \qquad \text{ if } m \text{ is even and } n \text{ is odd}. \end{cases}$$
Hence, we obtain
$$\sum_{k=0}^{m-1} \, \frac{p_{T, F}(m,n)}{m} = \begin{cases} \, p_{T, F}(m, 0) \qquad \qquad \qquad \qquad \text{if } m \text{ is odd},\\ \\ \, \displaystyle{\frac{p_{T, F}(m, 0)+p_{T, F}(m, 1)}{2} \qquad \text{if } m \text{ is even}.}\end{cases}$$
and (\ref{eq: 2.1}) follows.

Now, we consider the case of a reversal system and prove Theorem A.

\begin{prop}\label{prop: 3.1}
Let 
$$\mathcal{F}_1 = \{\langle a^m b^{2l}, b^{2k} \rangle : m>0, \, 0 \leq l<k \leq r \text{ and } k \, |\, r\}$$
and
$$\mathcal{F}_2 = \{\langle a^m, a^jb^{2k-1} \rangle : m>0, \, 0\leq j < m \text{ and } 2k-1 \, |\, r\}.$$
The set of finite index subgroups of $G_{2r}$ is a disjoint union of  $\mathcal{F}_1$ and $\mathcal{F}_2$.
\end{prop}

\begin{rmk}
If we set $H=\langle a, b^2\rangle$, then $H$ is an index 2 subgroup of $G_{2r}$. The set $\mathcal{F}_1$ is the collection of all the finite index subgroups of $H$. Thus, $K \in \mathcal{F}_1$ acts on $X$ by automorphisms. 
\end{rmk}

Before proving Proposition \ref{prop: 3.1}, we introduce some notation. For $i=1, \cdots, 2r-1$, We denote the order of $b^i$ by $|i|$:
$$|i|= \min\{n \in \mathbb{N} \setminus \{0\}: (b^i)^n = 1\},$$ 
or equivalently, $|i|= {2r}/{\gcd(i,2r)}$.

\begin{pf}
Since the proof is quite long, we will break it down into five steps and prove the following:
\newline
(1) If $K \in \mathcal{F}_1 \cup \mathcal{F}_2$, then $K$ is a finite index subgroup of $G_{2r}$.
\newline
(2) If $K=\langle a^mb^{2l}, b^{2k}\rangle$ is a finite index subgroup of $G_{2r}$ for some integers $m$, $l$ and $k$, then we may assume that $m>0$, $0 \leq l < k \leq r$ and $k \,| \, r$.
\newline
(3) If $K$ is a finite index subgroup of $\langle a, b^2 \rangle$, then $K \in \mathcal{F}_1$.
\newline
(4) If $K=\langle a^m, a^jb^{2k-1}\rangle$ is a finite index subgroup of $G_{2r}$, then we may assume that $m>0$, $0 \leq j < m$ and $2k-1 \, | \, r$.
\newline
(5) If $K$ is a finite index subgroup of $G_{2r}$ and $K\notin \mathcal{F}_1$, then $K \in \mathcal{F}_2$.

(1) Every finite index subgroup $K$ of $G$ must contain $a^m$ for some positive integer $m$. 
Suppose that $g\in G$ and that $a^m \in \langle g \rangle$. Then there is an integer $n$ such that $g^n=a^m$. We may assume that $n>0$ by replacing $g$ with $g^{-1}$ if $n<0$. 
It is clear that $m$ divides $n$. 
Then $g$ is either $g=a^{\frac{n}{m}}$ or $g=a^{\frac{n}{m}}b^{2l}$ for some positive integer $l$ with $|2l|\, \big| \, n$. 
Thus, if $K$ is generated by one element of $G$, then $K$ is either $\langle a^m \rangle$ for some positive integer $m$ or $\langle a^mb^{2l}\rangle$ for some positive integers $m$ and $l$. 
We see that the elements in $\mathcal{F}_1$ and $\mathcal{F}_2$ are finite index subgroups of $G_{2r}$.

(2) Clearly, we may assume that $m>0$, $0 \leq l, k \leq r$ and $k | r$. 
If $k=r$, then we are done. So, we assume that $0<k<r$.
If $k \, | \, l$, then $\langle a^mb^{2l}, b^{2k} \rangle = \langle a^m, b^{2k} \rangle \in \mathcal{F}_1$. 
To avoid duplication, we assume $k \nmid l$. By the division algorithm, there are integers $q$ and $p$ such that $l=qk+p$ and $0 < p < k$. 
Since $a^m b^{2p} = a^m b^{2l-2qk}$, it follows that
 $\langle a^mb^{2p}, b^{2k} \rangle = \langle a^m b^{2l}, b^{2k} \rangle$. Thus, it suffices to consider the case $0 \leq l<k \leq r$.

(3) We shall show that if $K$ is a finite index subgroup of $H=\langle a, b^2\rangle$, then $K$ belongs to $\mathcal{F}_1$. More precisely, we shall show that (i) $\langle a^m, a^n\rangle$ (ii) $\langle a^m, a^nb^{2l}\rangle$ and (iii) $\langle a^m b^{2k}, a^n b^{2l}\rangle$ belong to $\mathcal{F}_1$.

(i) B$\acute{\text{e}}$zout's identity tells us that for any positive integers $m$ and $n$, there are integers $k_1$ and $k_2$ such that $\gcd(m,n)=k_1 m+k_2 n$. Thus, if we set $d=\gcd(m, n)$, then we obtain $\langle a^m, a^n\rangle = \langle a^{d}\rangle$.

(ii) Suppose that $K=\langle a^m, a^nb^{2l}\rangle$ is not in $\mathcal{F}_1$  for some positive integers $m$, $n$ and $l$ with $0<l<r$.
If $m$ divides $n$, then $K = \langle a^m, b^{2l}\rangle \in \mathcal{F}_1$.
So we assume that $m$ does not divide $n$. 
If we set $m_1=\gcd(m, n|2l|)$, then $K=\langle a^{m_1}, a^nb^{2l}\rangle$ by B$\acute{\text{e}}$zout's identity. Again, if $m_1$ divides $n$, then $K \in \mathcal{F}_1$. Thus, $m_1$ does not divide $n$. 
We set 
$$m_{i+1} = \gcd(m_{i}, n|2l|) \qquad (i=1, 2, \cdots)$$
and continue this process.
Since $\{m_i\}$ is a finite decreasing sequence, there is a positive integer $d$ such that $m_d$ divides $n$ and we have $K=\langle a^{m_d}, b^{2l}\rangle \in \mathcal{F}_1$. This is a contradiction.
 
(iii) Suppose that $K=\langle a^m b^{2k}, a^n b^{2l}\rangle$ for some positive integers $m$, $n$, $k$ and $l$ with $0<k,\, l<r$. 
If we put $d=\gcd(m|2k|, n|2l|)$, then $a^d\in K$ by B$\acute{\text{e}}$zout's identity. 
It is clear that $K=\langle a^d, a^m b^{2k}, a^n b^{2l}\rangle$. 
By the argument in (ii), there are positive integers $i$ and $j$ such that $\langle a^d,a^mb^{2k} \rangle =\langle a^i, b^{2k}\rangle$ and $\langle a^d, a^nb^{2l}\rangle = \langle a^j, b^{2l}\rangle$.
Thus, we have $K=\langle a^i, a^j, b^{2k}, b^{2l}\rangle = \langle a^s, b^{2t}\rangle$ for some positive integers $s$ and $t$.
  
(4) We may assume that $2k-1$ is a divisor of $r$. 
If $j \geq m$, then there are positive integers $q$ and $p$ such that $j=qm+p$ and $0 \leq p<m$ by the division algorithm. 
Since $a^jb^{2k-1}=a^{qm+p}b^{2k-1}$, it follows that $K=\langle a^m, a^pb^{2k-1}\rangle$.

(5) We shall prove that (i) $\langle a^mb^{2l}, a^nb^{2k-1}\rangle$ and (ii) $\langle a^m, a^nb^{2k-1}, a^jb^{2l-1}\rangle$ belong to $\mathcal{F}_2$.

(i) Suppose that $K=\langle a^mb^{2l}, a^nb^{2k-1}\rangle$ and that $d=\gcd(2l, 2k-1)$. We shall show that $K=\langle a^m, a^nb^d\rangle$.
We note that $d$ is an odd number and that $2d$ divides $2l$. By B$\acute{\text{e}}$zout's identity, there are integers $k_1$ and $k_2$ such that $d=2l k_1 + (2k-1) k_2$. Since $d$ is odd, we see that $k_2$ must be odd. 

Since $(2k-1)/d$ and $d$ are odd, it follows that $a^nb^{2k-1}=(a^nb^d)^{(2k-1)/d}$ and that $b^{2d}=(a^n b^d)^2$.
Hence, $a^nb^{2k-1}, b^{2d} \in \langle a^m, a^nb^d\rangle$.
From $b^{2d}\in \langle a^m, a^nb^d\rangle$ and the fact that $2d$ divides $2l$, we obtain $b^{2l} \in \langle a^m, a^nb^d\rangle$ and 
we see that $a^m b^{2l} \in \langle a^m, a^nb^d\rangle$. 

Conversely, since $k_2$ is odd, it follows that
$$a^{mk_1 + n}b^d = a^{mk_1+n}b^{2lk_1 + (2k-1)k_2} = a^{mk_1} b^{2lk_1}a^nb^{(2k-1)k_2} = (a^m b^{2l})^{k_1}(a^n b^{2k-1})^{k_2}$$
and we have $a^{mk_1+n}b^d \in K$.
Since $d$ is odd, we have $b^{2d} = (a^{mk_1+n}b^d)^2 \in K$.
From this and the fact that $2d$ divides $2l$, we have $b^{2l} \in K$.
This implies $a^m \in K$.
Since $a^{mk_1+n}b^d, a^m \in K$, we have $\langle a^m, a^nb^d \rangle \subset K$.

(ii) Suppose that $K=\langle a^m, a^nb^{2k-1}, a^jb^{2l-1}\rangle$.
Since $a^{n-j}b^{2k+2l-2} \in K$ and $\langle a^nb^{2k-1}, a^jb^{2l-1}\rangle = \langle a^{n-j}b^{2k+2l-2}, a^nb^{2k-1}\rangle$, the result follows from (i).
$\hfill \Box$
\end{pf}

A direct calculation yields the following lemma.

\begin{lemma} \label{lemma: 3.2}
The indexes of the subgroups in $\mathcal{F}_1$ and $\mathcal{F}_2$ are as follows:
\newline
\emph{(1)} $|G/\langle a^mb^{2l}, b^{2k}\rangle|=2km \qquad (m>0, 0 \leq l<k \leq r \emph{ and } k \, |\, r)$.
\newline
\emph{(2)} $|G/\langle a^m, a^jb^{2k-1}\rangle|=(2k-1)m \qquad (m>0, 0\leq j <n \emph{ and } 2k-1 \, |\, r)$.
\end{lemma}

By Proposition \ref{prop: 3.1}, the Lind zeta function of $(X, T, R)$ can be decomposed into the product of
two terms:
$$\sum_{K\in \mathcal{F}_1} \frac{f_{\alpha}(K)}{|G_{2r}/K|} \,t^{|G_{2r}/K|} \qquad \text{and} \qquad \sum_{K\in \mathcal{F}_2} \frac{f_{\alpha}(K)}{|G_{2r}/K|} \,t^{|G_{2r}/K|}.$$

We introduce some notation.
We denote by $X_{2l}$ the set of points in $X$ fixed by $R^{2l}$:
$$X_{2l} = \{ x \in X: R^{2l}(x)=x \} \qquad (0 < l < r).$$ 
Dropping `$|_{X_{2l}}$', we denote the restrictions of $T$ and $R$ to $X_{2l}$ by $T$ and $R$, respectively.
For instance, if $(X, T, R)$ is a reversal system of order 12, 
it has 3 sub-reversal systems $(X_4, T, R)$, $(X_6, T, R)$ and $(X_{12}, T, R)=(X, T, R)$ and 2 sub-flip systems $(X_2, T, R)$ and $(X_6, T, R^3)$. 

If $m$ and $l$ are positive integers with $0 \leq l <r$, then
$f_{T, R}(m, 2l)$ will denote the number of points in $X$ fixed by $T^m \circ R^{2l}$:
$$f_{T, R}(m, 2l) = |\{x \in X : T^m \circ R^{2l}(x) = x\}|$$
and we set
\begin{equation}\label{eq: 3.1}
g_{T, R}(t) = \sum_{m=1}^{\infty}\, \sum_{l=0}^{r-1}\, \frac{f_{T, R}(m, 2l)}{m} \, t^{m}.
\end{equation}
Abusing notation, we refer to this function as a generating function. To avoid confusion, we denote it by $g$, while the generating function of a flip system is denoted by $h$.

If we denote by $g_{2k}$ the generating function of a sub-reversal system $(X_{2k}, T, R)$, we obtain
$$\sum_{K\in \mathcal{F}_1} \frac{f_{\alpha}(K)}{|G_{2r}/K|} \,t^{|G_{2r}/K|}=
\prod_{\substack{k|r \\ 1 \leq k \leq r}}\, \exp \Big(\frac{g_{2k}(t^{2k})}{2k} \Big).$$

If $x$ is fixed by $T^j \circ R^{2k-1}$, then $x \in X_{4k-2}$.
If we denote the generating function of a sub-flip system $(X_{4k-2}, T, R^{2k-1})$  by $h_{4k-2}$, then by the same argument in the sketch of the proof of (\ref{eq: 2.1}), we obtain
$$\sum_{K\in \mathcal{F}_2} \frac{f_{\alpha}(K)}{|G_{2r}/K|} \,t^{|G_{2r}/K|}=\prod_{\substack{2k-1|r \\ 1 \leq 2k-1 \leq r} }\exp \left( \frac{h_{4k-2}(t^{2k-1})}{2k-1} \right).$$
This completes the proof of Theorem A. 

\section{Proof of Theorem B}
A reversal $\varphi$ for a shift dynamical system $(X, \sigma_X)$ is said to be \textit{one-block} if
$$x, x' \in X \;\; \text{and} \;\; x_0=x'_0 \quad \Rightarrow \quad \varphi(x)_0=\varphi(x')_0.$$
When $\varphi$ is a one-block reversal of order $2r$, there is a unique map $\tau:\mathcal{B}_1(X) \rightarrow \mathcal{B}_1(X)$ such that $\tau ^{2r} = \text{id}_{\mathcal{B}_1(X)}$ and that
$$\varphi(x)_i=\tau(x_{-i}) \qquad (x \in X;\, i \in \mathbb{Z}).$$
We call $\tau$ the \textit{symbol map} of $\varphi$.

\begin{prop}\label{prop: 4.1}
Let $(X, \sigma_X, \varphi)$ be a shift-reversal system of order $2r$. Then there is a shift-reversal system $(Y, \sigma_Y, \psi)$ such that it is $G_{2r}$-conjugate to $(X, \sigma_X, \varphi)$ and
$\psi$ is one-block.
\end{prop}

\begin{pf}
Let $\mathcal{A}$ be the set of ordered $2r$-tuple $(a_{0} , a_{1}, \cdots, a_{2r-1})$ having the following property:
$$\exists x\in X \quad \text{s.t.} \quad a_{k}=\varphi^k(x)_0 \;\; \text{for all} \;\; k=0, 1, \cdots, 2r-1.$$
Obviously, $\mathcal{A}$ is finite.
We define a block map $\tau : \mathcal{A} \rightarrow \mathcal{A}$ and a map $\Psi : \mathcal{A}^{\mathbb{Z}} \rightarrow \mathcal{A}^{\mathbb{Z}}$ by
$$\tau(a_{0}, a_{1}, \cdots a_{2r-1})=(a_{1}, a_2, \cdots a_{2r-1}, a_{0}) \qquad ((a_{0}, a_{1}, \cdots a_{2r-1})\in \mathcal{A})$$
and
$$\Psi(y)_i = \tau(y_{-i}) \qquad (y \in \mathcal{A}^{\mathbb{Z}};\, i \in \mathbb{Z}).$$
so that $\tau^{2r} = \text{id}_{\mathcal{A}}$ and that $\Psi$ is a one-block reversal for $(\mathcal{A}^{\mathbb{Z}}, \sigma)$ of order $2r$.
Let $\theta : X \rightarrow \mathcal{A}^{\mathbb{Z}}$ be the map defined by
$\theta(x)_i=(a_0, a_1, \cdots, a_{2r-1})$ satisfying
$$a_k = \begin{cases} \varphi^k(x)_i \qquad \text{if } k \text{ is  even}, \\\varphi^k(x)_{-i} \qquad \text{if } k \text{ is  odd} \qquad(k=0, 1, \cdots, 2r). \end{cases}$$
Then $\theta$ is one-to-one and continuous and we have
$$\theta \circ {\sigma}_X = \sigma \circ \theta \qquad \text{and} \qquad \theta \circ \varphi = \Psi \circ \theta.$$
Setting $Y=\theta(X)$ and $\psi=\Psi|_Y$ yields the desired result.
\hfill $\Box$
\end{pf}

If $\mathcal{A}$ is a finite set and $A$ is a zero-one $\mathcal{A} \times \mathcal{A}$ matrix, then
$\textsf{X}_A$ will denote the topological Markov chain determined by $A$:
$$\textsf{X}_A=\{ x\in \mathcal{A}^{\mathbb{Z}} \,:\, \forall i\in \mathbb{Z} \;\; A(x_i, x_{i+1})=1 \}.$$
We denote the restriction of the shift map of $\mathcal{A}^{\mathbb{Z}}$ to $\textsf{X}_A$ by $\sigma_A$.
Suppose that $\mathcal{A}$ is a finite set and $A$ and $J$ are zero-one $\mathcal{A} \times \mathcal{A}$ matrices such that
\begin{equation}\label{eq: 4.1}
AJ=JA^{\textsf{T}} \qquad \text{and} \qquad J^{2r}=I.
\end{equation}
Since $J$ is a nonsingular zero-one matrix, it follows that there is a unique permutation $\tau_J$ of $\mathcal{A}$ such that
$$J(a , b) = 1  \qquad \Leftrightarrow \qquad \tau_J(a)=b \qquad (a, b \in \mathcal{A}).$$
From $J^{2r}=I,$ we have $(\tau_J)^{2r}=\text{id}_{\mathcal{A}}$. 
From $AJ=JA^{\textsf{T}}$, we have
\begin{equation}\label{eq: 4.2}
A(a,b)=A(\tau_J(b), \tau_J(a))\qquad(a, b \in \mathcal{A}).
\end{equation}
Now, we define $\varphi_J : \mathcal{A}^{\mathbb{Z}} \rightarrow \mathcal{A}^{\mathbb{Z}}$ by
$$\varphi_J(x)_i=\tau_J(x_{-i}) \qquad (x \in \mathcal{A}^{\mathbb{Z}};\, i \in \mathbb{Z}).$$
The map $\varphi_J$ is a one-block reversal for $(\mathcal{A}^{\mathbb{Z}}, \sigma)$ and has $\tau_J$ as a symbol map. 
Since (\ref{eq: 4.2}) implies that $\varphi_J(\textsf{X}_A)=\textsf{X}_A$,
the restriction $\varphi_{J, A}$ of $\varphi_J$ to $\textsf{X}_A$ is a reversal for $(\textsf{X}_A, \sigma_A)$.

The following is a corollary of Proposition \ref{prop: 4.1}.

\begin{cor}
Let $(X, \sigma_X, \varphi)$ be a shift-reversal system of finite type of order $2r$. Then there are zero-one square matrices $A$ and $J$ satisfying \emph{(\ref{eq: 4.1})} such that $(\textsf{X}_A, \sigma_A, \varphi_{J,A})$ is $G_{2r}$-conjugate to $(X, \sigma_X, \varphi)$.
\end{cor}

\begin{pf}
By Proposition \ref{prop: 4.1}, we may assume that $\varphi$ is a one-block reversal.
Since $(X, \sigma_X)$ is a shift of finite type, there is a positive number $N$ such that if $n>N$, then for all
$a_1 a_2 \cdots a_n$ and $b_1 b_2 \cdots b_n$ in $\mathcal{B}_n(Y)$, 
$$a_2 \cdots a_n = b_1 \cdots b_{n-1} \quad \text{implies} \quad a_1 a_2 \cdots a_n b_n =a_1 b_1 b_2 \cdots b_n\in \mathcal{B}_{n+1}(Y).$$
We fix an odd number $n>N$.
We define $\mathcal{B}_n(X) \times \mathcal{B}_n(X)$ matrices $A$ and $J$ by
$$A(a_1 \cdots a_n, \;b_1 \cdots b_n) = \begin{cases} 1 \qquad \text{if } a_2\cdots a_n = b_1 \cdots b_{n-1},\\ 0 \qquad \text{otherwise}\end{cases}$$
and
$$J(a_1 \cdots a_n, \;b_1 \cdots b_n) = \begin{cases} 1 \qquad \text{if } b_i=\tau(a_{n+1-i}) \quad \text{for all } i=1, 2, \cdots, n,\\ 0 \qquad \text{otherwise}.\end{cases}$$

Now one can show that the mapping
$$\cdots a_{-1} \underline{a_0} a_1 \cdots \qquad \mapsto \qquad \cdots \left[\begin{array}{c} a_{n-2} \\ \vdots \\ a_{-1} \end{array}\right] \underline{\left[\begin{array}{c} a_{n-1} \\ \vdots \\ a_{0} \end{array}\right]} \left[\begin{array}{c} a_{n} \\ \vdots \\ a_{1} \end{array}\right] \cdots$$
(Here, the zero-th coordinates are underlined.)
\newline
is a $G_{2r}$-conjugacy from
$(X, \sigma_X, \sigma_X^{n-1} \circ \varphi)$ to $(\textsf{X}_A, \sigma_A, \varphi_{J, A})$. 
Since $\sigma_X^{(n-1)/2}$ is a $G_{2r}$-conjugacy from $(X, \sigma_X, \varphi)$ to $(X, \sigma_X, \sigma_X^{n-1}\circ \psi)$ (See the third paragraph of Section 2), the result follows.
\hfill $\Box$
\end{pf}

In the rest of the section, we prove Theorem B. 
We first note that $A^{m}J^{2l}=J^{2l}A^m$ for all $m$ and $l$.
We shall prove 
$$f_{\sigma_A, \varphi_{J, A}}(m, 2l) =\text{tr}(J^{2l}A^m) \qquad (m=1,2, \cdots, 0 \leq l <r).$$
From the fact that
$$(\sigma_A)^m \circ (\varphi_{J, A})^{2l} (x)=x \quad \Leftrightarrow \quad x_{i}=(\tau_J)^{2l}(x_{m+i}) \qquad(x\in X; i \in \mathbb{Z})$$
and (\ref{eq: 4.2}), it follows that
$$f_{\sigma_A, \varphi_{J, A}}(m, 2l) 
= \big|\{x_0 \cdots x_{m-1} \in \mathcal{B}_m(\textsf{X}_A) : A \big((\tau_J)^{2l}(x_{m-1}), \; x_{0}\big)=1 \}\big|.$$
Since $\tau_J(a)$ is the unique element of $\mathcal{B}_1(\textsf{X}_A)$ satisfying $J(a, \tau_J(a))=1$ for all $a \in \mathcal{B}_1(\textsf{X}_A)$, 
it follows that
$$A\big((\tau_J)^{2l}(a),\; b \big)=J^{2l}A(a, b) \qquad(a, b \in \mathcal{A} ; \; 0 \leq l <r).$$
Hence 
$$f_{\sigma_A, \varphi_{J, A}}(m, 2l)
= \big|\{x_{m-1} \in \mathcal{B}_1(\textsf{X}_A) : J^{2l}A^m (x_{m-1}, x_{m-1})=1 \}\big|$$
and Theorem B is proved.

\section{Proof of Theorem C}

In this section, we assume that $(X, \sigma_X, \varphi)$ is a sofic-reversal system of order $2r$. By Proposition \ref{prop: 4.1}, we may assume that $\varphi$ is a one-block reversal. We will denote the symbol map of $\varphi$ by $\tau$.
We also assume that $\mathcal{A}$ is a finite set, $\mathcal{L}:\mathcal{A} \rightarrow \mathcal{B}_1(X)$ is a labeling and that $A$ and $J$ are zero-one $\mathcal{A} \times \mathcal{A}$ matrices having the following properties:
\newline
(P1) $A$ and $J$ satisfy (\ref{eq: 4.1}).
\newline
(P2) $\mathcal{L} \circ \tau_J = \tau \circ \mathcal{L}$.
\newline
(P3) $\mathcal{L}_{\infty}$ has no graph diamonds.
\newline
Recall that a labeling $\mathcal{L}:\mathcal{A} \rightarrow \mathcal{B}_1(X)$ gives a one-block code $\mathcal{L}_{\infty}$ from $\textsf{X}_A$ onto $X$, that is, $\mathcal{L}_{\infty}(y)_i=\mathcal{L}(y_i)$ for all $i \in \mathbb{Z}$. Obviously, $\mathcal{L}_{\infty}$ satisfies
\begin{equation}\label{eq: 5.1}
\mathcal{L}_{\infty} \circ \sigma_A = \sigma_X \circ \mathcal{L}_{\infty}.
\end{equation}
It is also obvious that (P2) implies  
\begin{equation}\label{eq: 5.2}
\mathcal{L}_{\infty} \circ \varphi_{J,A} = \varphi \circ \mathcal{L}_{\infty}.
\end{equation}
A graph diamond for $\mathcal{L}_{\infty}$ is a pair $(a_1 a_2 \cdots a_n ,$ $b_1 b_2 \cdots b_n)$ of distinct blocks in $\mathcal{B}(\textsf{X}_A)$
such that $a_1=b_1, a_n=b_n$ and $\mathcal{L}(a_i)=\mathcal{L}(b_i)$ for all $i = 1, \cdots, n$.

In \cite{KR}, it is proved that Kriger's joint state chain, found in \cite{K}, induces a presentation $(\mathcal{A}, \mathcal{L}, A, J)$ of $(X, \sigma_X, \varphi)$ which has the properties (P1)-(P3) when $\varphi$ is a one-block flip. We will show that this is still valid when $\varphi$ is a one-block reversal at the end of the section.

Our proof of Theorem C is basically a slight modification of the proof of (\ref{eq: 2.3}) described in Section 6.4. of \cite{LM}. 
In \cite{KR}, the authors proved (\ref{eq: 2.3}) as well to modify it and prove (\ref{eq: 2.5}).
We adopt some notation and arguments introduced in \cite{KR}.

We will denote by $P_i$ the projection of $\textsf{X}_A$ onto $\mathcal{A}$ which assigns to each point in $\textsf{X}_A$ the $i$-th coordinate of it:
$$y \mapsto y_i \qquad (y \in \textsf{X}_A)$$
and we will denote the cardinality of a set $S$ by $|S|$.

The following is an immediate consequence of Property (P3).

\begin{lemma}\label{lem: 5.1}
If $x$ is periodic, then $P_i$ is one-to-one on $\mathcal{L}_{\infty}^{-1}(x)$ for every $i$. 
As a result, if $x$ is periodic, then $1 \leq |\mathcal{L}_{\infty}^{-1}(x)| \leq |\mathcal{A}|$.
\end{lemma}

Let $F(m, 2l)$ be the set of points in $X$ fixed by $(\sigma_X)^m \circ \varphi^{2l}$:
$$F(m , 2l) = \{ x \in X : ((\sigma_X)^m \circ \varphi^{2l})(x) = x \}$$
so that $f_{\sigma_X, \varphi}(m, 2l) = |F(m, 2l)|$.
If $x \in F(m, 2l)$, then $x$ is a periodic point of period $m|2l|$ and we have $1 \leq |\mathcal{L}_{\infty}^{-1}(x)| \leq |\mathcal{A}|$ by Lemma \ref{lem: 5.1}.
From (\ref{eq: 5.1}) and (\ref{eq: 5.2}), it follows that $((\sigma_A)^m\circ (\varphi_{J, A})^{2l}) (\mathcal{L}_{\infty}^{-1}(x))=\mathcal{L}_{\infty}^{-1}(x)$.
This implies that the restriction of $(\sigma_A)^m \circ (\varphi_{J, A})^{2l}$ to $\mathcal{L}_{\infty}^{-1}(x)$ is a permutation of $\mathcal{L}_{\infty}^{-1}(x)$.
When $\pi$ is a permutation of a finite set, its sign will be denoted by sgn$(\pi)$. 

The following is proved in Section 6.4 of \cite{LM}.
\begin{lemma} \label{lemma: 5.2}
Let $\pi$ be a permutation of a non-empty finite set $S$ and $\mathcal{C} = \{ E \subset S : \pi(E) = E \}.$
Then 
$$\sum_{E \in \mathcal{C}\setminus \{\varnothing\}} (-1)^{|E|+1} \, \emph{sgn} (\pi |_E)=1.$$
\end{lemma}

For $x \in F(m, 2l)$, we denote by $\mathcal{C}(m, 2l; \,x)$ the set of all $(\sigma_A)^m \circ (\varphi_{J, A})^{2l}$-invariant subsets of $\mathcal{L}_{\infty}^{-1}(x)$:
$$\mathcal{C}(m, 2l;\, x)=\{E \subset \mathcal{L}_{\infty}^{-1}(x) : ((\sigma_A)^m\circ {\varphi}^{2l}) (E) =E \}.$$ 
Now, Lemma \ref{lemma: 5.2} tells us that
$$f_{\sigma_X, \varphi}(m, 2l) = \sum_{x \in F(m, 2l)} \sum_{E \in \mathcal{C}(m, 2l;\, x)\setminus \{ \varnothing \} } (-1)^{|E|+1} \text{sgn}((\sigma_A)^m \circ (\varphi_{J, A})^{2l}|_E).$$
If $x$ and $x'$ are distinct points in $F(m, 2l)$, then $\mathcal{C}(m, 2l;\, x) \cap \mathcal{C}(m, 2l;\, x') = \{ \varnothing \}$. 
Thus, we set 
$$\mathcal{C}(m, 2l)=\bigcup_{x \in F(m, 2l)} \mathcal{C}(m, 2l;\, x)$$ 
and write 
$$f_{\sigma_X, \varphi}(m, 2l) = \sum _{E \in \mathcal{C}(m, 2l) \setminus \{ \varnothing \}} (-1)^{|E|+1} \text{sgn}((\sigma_A)^m\circ(\varphi_{J, A})^{2l}|_E).$$
On the other hand, if we set
$$\mathcal{C}_k(m, 2l) = \{ E \in \mathcal{C}(m, 2l) : |E| = k \} \qquad (k=1, 2, \cdots, |\mathcal{A}|),$$
then $\mathcal{C}(m, 2l) \setminus \{\varnothing \}$ is a disjoint union of
$\mathcal{C}_1(m, 2l)$, $\mathcal{C}_2(m, 2l)$, $\cdots$, $\mathcal{C}_{|\mathcal{A}|}(m, 2l)$
and we obtain
$$f_{\sigma_X, \varphi}(m, 2l)=\sum_{k=1}^{|\mathcal{A}|} (-1)^{k+1} \sum_{E \in \mathcal{C}_k(m, 2l)} \text{sgn}((\sigma_A)^m\circ (\varphi_{J, A})^{2l}|_E).$$

Now, we shall construct matrices $A_k$ and $J_k$ satisfying 
\begin{equation}\label{eq: 5.3}
\sum_{E \in \mathcal{C}_k(m, 2l)} \text{sgn}((\sigma_A)^m\circ (\varphi_{J, A})^{2l}|_E)=\text{tr}((A_k)^{m}(J_k)^{2l})
\end{equation}
for all $k=1, \cdots |\mathcal{A}|$.

We fix a linear order $<$ of $\mathcal{A}$.
If $f$ is a one-to-one function from a subset $S$ of $\mathcal{A}$ into $\mathcal{A}$, we define 
$$N(<, \,f) = |\{ (a, b) \in S \times S : a < b \quad \text{and} \quad f(b) < f(a) \}|$$
and
$$\text{sgn}_<(f) = (-1)^{N(<, \, f)}.$$
If $g : f(S) \rightarrow \mathcal{A}$ is one-to-one again, then
$$\text{sgn}_<(g \circ f) = \text{sgn}_<(g)\, \text{sgn}_<(f).$$
From this, we see that 
\begin{equation}\label{eq: 5.11}
\text{sgn}_<(f^{-1}) = \text{sgn}_<(f).
\end{equation}
In particular, if $\pi$ is a permutation of a subset of $\mathcal{A}$, its sign is equal to $\text{sgn}_<(\pi)$.
From now on, we drop $<$ and write $\text{sgn}(f) = \text{sgn}_<(f)$.

For each positive integer $k \leq |\mathcal{A}|$, we set
$$\mathcal{A}_k = \{ S \subset \mathcal{A} : |S| = k \quad \text{and} \quad |\mathcal{L}(S)|=1 \}.$$ 
For $S_1, S_2 \in \mathcal{A}_k$, we denote by $F(S_1, S_2)$ the set of one-to-one correspondence $f : S_1 \rightarrow S_2$ such that
$A(a, f(a))=1$ for all $a \in S_1$. 

Define $\mathcal{A}_k \times \mathcal{A}_k$ matrix $A_k$ and $J_k$ by
$$A_k(S_1, S_2) = \sum_{f \in F(S_1, S_2)} \text{sgn}(f).$$
$$J_k(S_1, S_2) = \begin{cases} \text{sgn}(\tau_J|_{S_1}) \qquad \text{if } \tau_J(S_1)=S_2, \\ 0 \qquad \qquad \qquad \text{otherwise}. \end{cases}$$
It is straightforward to see that $(J_k)^{2r}=I$. We shall show that 
\begin{equation}\label{eq: 5.13}
A_kJ_k=J_kA_k^{\textsf{T}}.
\end{equation}
When $S, T \subset \mathcal{A}$ and $f: S \rightarrow T$ is a one-to-one correspondence, we define the map
$\tau_J(f): \tau_J(T) \rightarrow \tau_J(S)$ by
$$\tau_J(a) \mapsto \tau_J(f^{-1}(a)) \qquad (a \in T)$$
so that 
\begin{equation}\label{eq: 5.4}
\tau_J(f) = (\tau_J|_S) \circ f^{-1} \circ({\tau_J}^{2r-1}|_{\tau_J(T)}).
\end{equation}
It is clear that $\tau_J(f)$ is a one-to-one correspondence and that $(\tau_J)^{2r}(f)=f$.
From (\ref{eq: 5.4}), it follows that
\begin{equation}\label{eq: 5.12}
\text{sgn}(\tau_J(f)) = \text{sgn}(\tau_J|_S)\, \text{sgn}(f^{-1})\,\text{sgn}((\tau_J)^{2r-1}|_{\tau_J(T)}).
\end{equation}
By Property (P2), $S \in \mathcal{A}_k$ implies $\tau_J(S) \in \mathcal{A}_k$
and by (\ref{eq: 4.2}), $f \in F(S_1, S_2)$ implies that $\tau_J(f) 
\in F(\tau_J(S_2), \tau_J(S_1))$.
It is clear that the map $f \mapsto \tau_J(f)$ is a one-to-one correspondence.
From (\ref{eq: 5.11}), (\ref{eq: 5.12}) and the fact that $f \mapsto \tau_J(f)$ is a one-to-one correspondence, we obtain
$$A_k(\tau_J(T), \tau_J(S)) = (J_k)^{2r-1}(\tau_J(S), S)\,A_k(S, T)\,J_k(T, \tau_J(T)) \qquad(S, T \in \mathcal{A}_k)$$
and this proves (\ref{eq: 5.13}).

The followgin is an immediate consequence of Property (P3). It is also proved in \cite{KR}.
\begin{lemma}\label{lemma: 5.3}
Suppose that $S_0, S_1, \cdots, S_m \in \mathcal{A}_k$ and that $S_0=S_m$.
Then the following are equivalent.
\newline
\emph{(a)} $\prod_{j=0}^{m-1}A_k(S_j, S_{j+1}) \neq 0$.
\newline
\emph{(b)} $|F(S_j, S_{j+1})| \geq 1$ for $0 \leq j \leq m-1$.
\newline
\emph{(c)} $|F(S_j, S_{j+1})| = 1$ for $0 \leq j \leq m-1$.
\end{lemma}

The following lemma is a modification of Lemma 2.6 in \cite{KR} and 
(\ref{eq: 5.3}) is proved by this lemma:

\begin{lemma}\label{lemma: 5.4}
If $E \in \mathcal{C}_k(m, 2l)$ and
\begin{equation}\label{eq: 5.5}
S_j=P_j(E) \qquad (j=0, \cdots, m),
\end{equation}
then $A_k(S_j, S_{j+1}) \in \{ -1, 1 \}$ for all $j=0, \cdots, m$, 
$S_0=(\tau_J)^{2l}(S_m)$ and
$$\emph{sgn}((\sigma_A)^m\circ (\varphi_{J, A})^{2l}|_E) = (A_k)^{m}(J_k)^{2l}(S_0, S_0).$$
Conversely, suppose that $S_0, S_1, \cdots, S_m \in \mathcal{A}_k$. If $A_k(S_j, S_{j+1}) \in \{ -1, 1 \}$ for all $j=0, \cdots, m$ and $S_0=(\tau_J)^{2l}(S_m)$, then there is a unique $E \in \mathcal{C}_k(m, 2l)$
such that \emph{(\ref{eq: 5.5})} holds.
\end{lemma}

\begin{pf}
Suppose that $E \in \mathcal{C}_k(m, 2l)$ and $S_j=P_j(E)$ for all $j \in \mathbb{Z}$. 
Then $E \subset \mathcal{L}_{\infty}^{-1}(x)$ for some $x \in F(m, 2l)$ 
and we have
\begin{equation}\label{eq: 5.6}
\mathcal{L}(S_j) = \{ \mathcal{L}(y_i) : y \in E \} = \{ x_j \} \qquad (j \in \mathbb{Z}).
\end{equation}
Lemma \ref{lem: 5.1} implies that $P_j|_E : E \rightarrow S_j$ is a one-to-one correspondence for all $j$.
Hence, we see that
$$|S_j|=|E|=k \qquad \text{and} \qquad S_j \in \mathcal{A}_k \qquad (j \in \mathbb{Z}).$$
Since $({\sigma_A})^m\circ (\varphi_{J, A})^{2l}(E)=E, $ we have $S_j=(\tau_J)^{2l}(S_{m+j})$ for all $j$.
If we set $f_j = P_{j+1} \circ (P_j|_E)^{-1}$, then $f_j \in F(S_j, S_{j+1})$ for all $j$.
Since $S_0 = S_{m|2l|}$, Lemma \ref{lemma: 5.3} implies that 
$F(S_j, S_{j+1}) = \{ f_j \}$ for all $j$
and this implies that $A_k(S_j, S_{j+1}) \in \{ -1, 1 \}$ for all $j$.
As a consequence, we have
\begin{eqnarray*}
\prod_{j=0}^{m-1} A_k (S_{j}, S_{j+1})(J_k)^{2l}(S_{m}, S_0) 
&=& \text{sgn} ((\tau_J)^{2l} \circ f_{m-1} \circ \cdots f_1 \circ f_{0})\\ 
&=& \text{sgn} ((\tau_J)^{2l} \circ P_{m} \circ (P_{0}|_E)^{-1}). 
\end{eqnarray*}
On the other hand, for all $y \in \textsf{X}_A$, we have $P_{m}(y)=P_{0}((\sigma_A)^m(y))$ and $(P_m \circ (\varphi_{J, A})^{2l})(y) = ((\tau_J)^{2l}\circ P_m)(y)$. These imply
\begin{eqnarray*}
\text{sgn}((\sigma_A)^m\circ(\varphi_{J, A})^{2l}|_E) 
&=& \text{sgn} (P_{0} \circ (\sigma_A)^m \circ (\varphi_{J, A})^{2l} \circ (P_{0}|_E)^{-1})\\
& =& \text{sgn} ((\tau_{J})^{2l}\circ P_{m} \circ (P_{0}|_E)^{-1}).
\end{eqnarray*}
This proves the first assertion.

To prove the second assertion, we assume that $S_j \in \mathcal{A}_k$ satisfy $A_k(S_j, S_{j+1}) \in \{ -1, 1 \}$ for all $j=0, \cdots, m$ and $S_0=(\tau_J)^{2l}(S_m)$. 
By (\ref{eq: 5.13}), $A_k(S_j, S_{j+1}) \in \{ -1, 1 \}$ if and only if $$A_k((\tau_J)^{2i}(S_j), (\tau_J)^{2i}(S_{j+1})) \in \{ -1, 1 \} \qquad (\forall \; i \in \mathbb{Z};\, j=0, \cdots, m).$$ 
We set for $j=0, \cdots, m-1$,
$$S_{m+j}=(\tau_J)^{|2l|-2l}(S_j),$$
$$S_{2m+j}=(\tau_J)^{|2l|-4l}(S_j),$$
$$\vdots$$
$$S_{(|2l|-2)m+j}=(\tau_J)^{4l}(S_j),$$
$$S_{(|2l|-1)m+j}=(\tau_J)^{2l}(S_j).$$
Since $A_k(S_{m-1}, S_m) \in \{ -1, 1 \}$ and $S_0=(\tau_J)^{2l}(S_m)$, it follows that 
$$A_k(S_j, S_{j+1}) \in \{ -1, 1 \} \qquad (j=0, 1, \cdots, m|2l|-1)$$ and 
$S_0=S_{m|2l|}$. 
Lemma \ref{lemma: 5.3} tells us that there are functions $f_0, f_1, \cdots, f_{m|2l|-1}$ such that $F(S_j, S_{j+1}) = \{ f_j \} $ for $j=0, \cdots, m|2l|-1$.
If $E$ denotes the set of all $y \in \textsf{X}_A$ such that
$$i \in \mathbb{Z} \;\;\text{and} \;\; 0 \leq j \leq m|2l|-1 \;\, \Rightarrow \;\, y_{m|2l|i+j} \in S_j \;\; \text{and} \;\; y_{m|2l|i+j+1}=f_j(y_{m|2l|i+j}),$$
then $E \in \mathcal{C}_k(m, 2l)$ and (\ref{eq: 5.5}) follows.

It remains to prove the uniqueness of $E$. 
Suppose that $x, x' \in F(m, 2l)$, $E, E' \in \mathcal{C}_k(m, 2l)$, $E \subset \mathcal{L}_{\infty}^{-1}(x)$, $E' \subset \mathcal{L}_{\infty}^{-1}(x')$
and that $P_j(E)=P_j(E')$ for $0 \leq j \leq m-1$. Then (\ref{eq: 5.6}) implies that $x=x'.$
Since $x$ and $x'$ are periodic points of period $m|2l|$, from Lemma \ref{lem: 5.1}, it follows that $E=E'$. This completes the proof.
\hfill $\Box$
\end{pf}

In the rest of the section, we show that Krieger's joint state chain of $(X, \sigma_X)$ has a natural one-block reversal and the resulting shift-reversal system has properties (P1)-(P3).

Let $X^+$ and $X^-$ denote the set of right-infinite sequences and left-infinite sequences which appear in $X$, respectively:
$$X^+ = \{ x_{[0, \infty)}: x \in X \} \qquad \text{and} \qquad X^- = \{ x_{(-\infty, 0]}: x \in X \}.$$
A \textit{future} is a subset $F$ of $X^+$ satisfying
$$\exists \lambda \in X^- \qquad \text{s.t.} \qquad F = \{ \rho \in X^+ : \lambda \rho \in X \}$$
and a \textit{past} is a subset $P$ of $X^-$ satisfying
$$\exists \rho \in X^+ \qquad \text{s.t.} \qquad P = \{ \lambda \in X^- : \lambda \rho \in X \}.$$
When $F$ is a future, $P$ is a past and $a \in \mathcal{B}_1(X)$,
we define $F(a)$ and $P(a)$ to be
$$F(a) = \{ \rho \in X^+ : a \rho \in F \} \qquad \text{and} \qquad P(a) = \{ \lambda \in X^- : \lambda a \in P \}.$$
A \textit{joint state} is a triple $(F, a, P)$ of a future $F$, a symbol $a$ in $X$ and a past $P$ with $F(a) \neq \varnothing$ and $P(a) \neq \varnothing$.
Let $\mathcal{A}$ denote the set of all joint states. 
It is well known \cite{K, LM} that there are finitely many futures and pasts and this implies that $\mathcal{A}$ is finite.
We define the labeling $\mathcal{L} : \mathcal{A} \rightarrow \mathcal{B}_1(X)$ by
$$\mathcal{L}(F, a, P) = a$$
and define the zero-one $\mathcal{A} \times \mathcal{A}$ matrix $A$ to be 
$$A((F_1, a_1, P_1), (F_2, a_2, P_2)) = \begin{cases} 1 \qquad F_1(a_1) = F_2 \quad \text{and} \quad P_1 = P_2(a_2), \\ 0 \qquad \text{otherwise}.\end{cases}$$
The topological Markov chain $\textsf{X}_A$ is called \textit{Krieger's joint state chain} of $X$.
By definition, $\mathcal{L}_{\infty}$ has no graph diamonds.

We define $\tau_{\infty}^+ : X^+ \rightarrow X^-$ by 
$$\rho = \rho_0 \, \rho_1 \, \rho_2 \cdots \mapsto \cdots \tau(\rho_2) \tau(\rho_1) \tau(\rho_0) \qquad (\rho \in X^+)$$
and define $\tau_{\infty}^- : X^- \rightarrow X^+$ by
$$\lambda = \cdots \lambda_{-2} \lambda_{-1} \lambda_0 \cdots \mapsto  \tau(\lambda_0) \tau(\lambda_{-1}) \tau(\lambda_{-2}) \cdots \qquad (\lambda \in X^-).$$
Since $\varphi(X)=X$, these maps are well-defined and we have $(F, a, P)$ is a joint state chain if and only if so is $(\tau_{\infty}^-(P), \tau(a), \tau_{\infty}^+(F))$. 
Since $\tau^{2r}$ is the identity map of $\mathcal{B}_1(X)$, it follows that $(\tau_{\infty}^- \circ \tau_{\infty^+})^r = \text{id}_{X^+}$ and $(\tau_{\infty}^+ \circ \tau_{\infty^-})^r = \text{id}_{X^-}$.
We define the zero-one $\mathcal{A} \times \mathcal{A}$ matrix $J$ to be
$$J((F_1, a_1, P_1), (F_2, a_2, P_2))= \begin{cases} 1 \qquad 
(\tau_{\infty}^-(P_1), \, \tau(a_1), \, \tau_{\infty}^+(F_1) )=(F_2, a_2, P_2),
\\0 \qquad \text{otherwise}.\end{cases}$$
It is straightforward to see that $(\mathcal{A}, \mathcal{L}, A, J)$ has properties (P1) and (P2).

\section{Some Remarks and Example}

In \cite{KR}, it is proved that the generating function $h_{\sigma_X, \varphi}(t)$ in  (\ref{eq: 2.11}) of a sofic shift-flip system $(X, \sigma_X, \varphi)$ is $\mathbb{N}$-rational. For the definition of $\mathbb{N}$-rationality, see \cite{R, BR1, E}.
We show that the generating function $g_{\sigma_X, \varphi}(t)$ in (\ref{eq: 3.1}) of a sofic shift-reversal system $(X, \sigma_X, \varphi)$ is $\mathbb{N}$-rational.
Since the set 
$$\bigcup_{m=1}^{\infty} \{w \in \mathcal{B}_m(X) : x_{[0, m-1]} = w \text{ for some } x \in F(m, 2l)\}$$
is recognizable, it follows that
$$\sum_{m=1}^{\infty} f_{\sigma_X, \varphi}(m, 2l) t^{m}$$
is $\mathbb{N}$-rational. 
Since $x\in F(m, 2l)$ if and only if $(\sigma_X)^i(x) \in F(m, 2l)$ for all
$i=1, \cdots, m-1$,
Berstel's theorem \cite{BR2, E, S} and Soitolla's theorem \cite{BR2, S} imply $\mathbb{N}$-rationality of $g_{\sigma_X, \varphi}(t)$.
The following example shows that neither $\frac{g_{2k}(t^{2k})}{2k}$ nor $\frac{h_{4k-2}(t^{2k-1})}{2k-1}$
is $\mathbb{N}$-rational in general. 

\begin{exa}
Let $\mathcal{A}=\{1, 2, \cdots, 7\}$ and let $A$ and $J$ be zero-one $\mathcal{A}\times\mathcal{A}$ matrices given by
$$A=\left[\begin{array}{rrrrrrr} 0& 1 & 0 & 0 & 0 & 1 & 1 \\ 0 & 0 & 0 & 0 & 0 & 0 & 1 \\ 0 & 1 & 0 &  1 & 0 & 0 & 1 \\ 0 & 0 & 0 & 0 & 0 & 0 & 1 \\ 0 & 0 & 0 & 1 & 0 & 1 & 1 \\ 0 & 0 & 0 & 0 & 0 & 0 & 1 \\ 1 & 1 & 1 & 1 & 1 & 1 & 1 \end{array}\right] \quad \text{and} \quad J=\left[\begin{array}{rrrrrrr} 0 & 1 & 0 & 0 & 0 & 0 & 0  \\ 0 & 0 & 1 & 0 & 0 & 0 & 0 \\ 0 & 0 & 0 & 1 & 0 & 0 & 0 \\ 0 & 0 & 0 & 0 & 1 & 0 & 0 \\ 0 & 0 & 0 & 0 & 0 & 1 & 0 \\ 1 & 0 & 0 & 0 & 0 & 0 & 0 \\ 0 & 0 &  0 & 0 & 0 & 0 & 1 \end{array}\right].$$
Then $(\textsf{X}_A, \sigma_A, \varphi_{J, A})$ is a shift-reversal system of finite type of order $6$.

We set 
$$\textsf{X}_{2l} = \{ x \in \textsf{X}_A: (\varphi_{J,A})^{2l}(x)=x \} \qquad (l=1,2,3).$$
Dropping `$|_{\textsf{X}_{2l}}$', we denote the restrictions of $\sigma_A$ and $\varphi_{J, A}$ to $\textsf{X}_{2l}$ by $\sigma_A$ and $\varphi_{J, A}$, respectively. 
Then $(\textsf{X}_A, \sigma_A, \varphi_{J, A})$ has one sub-reversal system $(\textsf{X}_6, \sigma_A, \varphi_{J, A})$ which is equal to itself and two sub-flip systems $(\textsf{X}_2, \sigma_A, \varphi_{J, A})$ and $(\textsf{X}_6, \sigma_A, (\varphi_{J, A})^3)$.  

It is easy to see that $$\textsf{X}_2 = \{7^{\infty}\}.$$
Hence, $$g_2(t) =\sum_{m=1}^{\infty} \,\frac{1}{m}\, t^m = \frac{1}{1-t}.$$
Direct computations show that
$$\text{tr}(A^mJ^2) = \text{tr}(A^mJ^4) = 1$$
and that
$$g_6(t) =\sum_{m=1}^{\infty} \,\frac{\text{tr}(A^m)}{m}\, t^m + \sum_{m=1}^{\infty} \,\frac{2}{m}\, t^m = \frac{1}{1-t-6t^2-6t^3}+\frac{2}{1-t}.$$

On the other hand, direct computations yield
$$h_6(t) = \frac{t + t^2 + 3 t^4 +3 t^6}{1 - t^2 - 6 t^4 -6 t^6} \qquad \text{and} \qquad h_2(t)=\frac{t}{1-t}.$$  

We note that $g_6(t^6)/6$, $g_2(t^2)/2$ and $h_6(t^3)/3$ are not $\mathbb{N}$-rational.
\end{exa}

\end{document}